\documentclass[12pt,reqno]{amsart}
\usepackage{enumerate, latexsym, amsmath, amsfonts, amssymb, amsthm, color}
\def\pmod #1{\ ({\rm{mod}}\ #1)}

\def\bg{\bigg}
\def\({\bg(}
\def\){\bg)}

\def\f{\frac}

\def\ls{\leqslant}

\def\bi{\binom}

\def\eq{\equiv}

\theoremstyle{plain}
\newtheorem{theorem}{Theorem}

\newtheorem{lemma}{Lemma}

\theoremstyle{definition}

\theoremstyle{remark}
\newtheorem{remark}{Remark}

 \vspace{4mm}

\begin{document}

\hbox{Preprint}
\medskip

\title
[{Proof of two supercongruences conjectured by Z.-W.Sun }]
{Proof of two supercongruences conjectured by Z.-W.Sun involving Catalan-Larcombe-French numbers}

\author
[Guo-Shuai Mao] {Guo-Shuai Mao}

\address {(Guo-Shuai Mao) Department of Mathematics, Nanjing
University, Nanjing 210093, People's Republic of China}
\email{mg1421007@smail.nju.edu.cn}
 \begin{abstract} The harmonic numbers $H_n=\sum_{0<k\ls n}1/k\ (n=0,1,2,\ldots)$ play important roles in mathematics. With helps of some combinatorial identities, we establish the following two congruences:
 $$\sum_{k=0}^{\frac{p-3}2}\f{\binom{2k}k^2H_k}{(2k+1)16^k}\  \mbox{modulo}\ p^2\ \mbox{and}\ \sum_{k=0}^{\frac{p-3}2}\f{\binom{2k}k^2H_{2k}}{(2k+1)16^k}\ \mbox{modulo}\ p$$
for any prime $p>3$, the second one was conjectured by Z.-W. Sun in 2012. These two congruences are very important to prove the following
 conjectures of Z.W.Sun: For any old prime $p$, we have
  $$\sum_{k=0}^{p-1}\frac{P_k}{8^k}\equiv1+2\big(\frac{-1}p\big)p^2E_{p-3}\pmod{p^3}$$
    and $$\sum_{k=0}^{p-1}\frac{P_k}{{16}^k}\equiv\big(\frac{-1}p\big)-p^2E_{p-3}\pmod{p^3},$$
  where $P_n=\sum_{k=0}^n\frac{\binom{2k}k^2\binom{2(n-k)}{n-k}^2}{\binom nk}$ is the n-th Catalan-Larcombe-French number.
\end{abstract}
\maketitle

\section{Introduction}
\setcounter{lemma}{0}
\setcounter{theorem}{0}
\setcounter{corollary}{0}
\setcounter{remark}{0}
\setcounter{equation}{0}
\setcounter{conjecture}{0}

The Catalan-Larcombe-French numbers $P_0,P_1,P_2,\ldots$(cf. \cite{FH}) are given by $$P_n=\sum_{k=0}^n\frac{\binom{2k}k^2\binom{2(n-k)}{n-k}^2}{\binom nk}=2^n\sum_{k=0}^{\lfloor n/2\rfloor}\binom{n}{2k}\binom{2k}k^24^{n-2k},$$ they arose from the theory of elliptic integrals (see \cite{PD}). It is known that $(n+1)P_{n+1}=(24n(n+1)+8)P_n-128n^2P_{n-1}$ for all $n\in \mathbb{Z}^{+}$. The sequence $(P_n)_{n\geq0}$ is also related to the theory of modular forms, see D.Zagier\cite{D}.

 Let $S_n=\frac{P_n}{2^n}$, Zagier noted that $$S_n=\sum_{k=0}^{[n/2]}\binom{2k}k^2\binom n{2k}4^{n-2k}.$$Recently Z.-W.Sun stated that
 \begin{align*}
 &S_n=\sum_{k=0}^n\binom{2k}k^2\binom k{n-k}(-4)^{n-k}
 \\=&\frac1{(-2)^n}\sum_{k=0}^n\binom{2k}k\binom{2(n-k)}{n-k}\binom k{n-k}(-4)^k.
 \end{align*}
  The Bernoulli numbers $\{B_n\}$  are defined by
$$B_0=1,\ \sum_{k=0}^{n-1}\binom nkB_k=0(n\geq2).$$
The Euler numbers$\{E_n\}$ are defined by
$$ E_0=1,E_n=-\sum_{k=1}^{[n/2]}\binom{n}{2k}E_{n-2k}(n\geq1),$$
where $[a]$ is the greatest integer not exceeding $a$.

{\bf Conjecture 1}\cite[Remark 3.13]{Su1} We also have the following conjecture related to Euler numbers:

$$\sum_{k=0}^{p-1}\frac{P_k}{8^k}\equiv1+2\big(\frac{-1}p\big)p^2E_{p-3}\pmod{p^3}$$
and
$$\sum_{k=0}^{p-1}\frac{P_k}{16^k}\equiv\big(\frac{-1}p\big)-p^2E_{p-3}\pmod{p^3}$$
 for any odd prime $p$, where $\big(\frac{\cdot}p\big)$ is the Legendre symbol. In this paper we mainly prove Conjecture 1.

 Recently Z.H. Sun and X.J Ji done some research on the Catalan-Larcombe-French numbers. In \cite{S2} Z.H. Sun showed that
 $$\sum_{k=0}^{p-1}\binom{2k}k\frac{S_k}{m^k}\pmod p\ \mbox{for}\  m =7,16,25,32,64,160,800,1600,156832.$$
 X.J Ji and Z.H. Sun \cite{JS} obtained $$S_{mp^r}\pmod{p^{r+2}},\ S_{mp^r-1}\pmod{p^r}\ \mbox{and}\  S_{mp^r+1}\pmod {p^{2r}},$$ where $p$ is an odd prime and $m, r$ are positive integers.

  In \cite{Su2} Z.-W.Sun showed that
  \begin{align}\label{1.1}
  \sum_{k=0}^{\frac{p-3}2}\frac{\binom{2k}k^2}{(2k+1)16^k}\equiv-2q_p(2)-pq_p(2)^2+\frac5{12}p^2B_{p-3}\pmod{p^3},
  \end{align}
  where $q_p(2)$ denotes the Fermat quotient $(2^{p-1}-1)/p$.

 Z-W.Sun \cite{Su3} also done some research on Harmonic numbers $H_n$ and $H_{2n}$, such as
 $$\sum_{k=0}^{(p-1)/2}\binom{2k}k^2\frac{H_k}{16^k}\equiv2\big(\frac{-1}p\big)H_{(p-1)/2}\pmod{p^2}$$
 and
 $$\sum_{k=0}^{(p-1)/2}\binom{2k}k^2\frac{H_{2k}}{16^k}\equiv\frac32\big(\frac{-1}p\big)H_{(p-1)/2}+pE_{p-3}\pmod{p^2}.$$

   Z.-W. Sun proved that
  \begin{align}\label{1.2}
  \sum_{k=0}^{(p-1)/2}\frac{\binom{2k}k^2}{16^k}\equiv(-1)^{(p-1)/2}+p^2E_{p-3}\pmod{p^3}
  \end{align}
  in \cite{Su4}, which is very important to get our results.

  Motivated by the above work, we mainly obtain the following result in this paper.
  \begin{theorem}\label{Th1.1} Let $p>3$ be a prime. Then
 \begin{equation}\label{1.3}
 \sum_{k=0}^{\frac{p-3}2}\f{\bi{2k}k^2H_k}{(2k+1)16^k}\eq4q_p(2)^2+2\big(\frac{-1}p\big)(E_{2p-4}-2E_{p-3})+\frac7{12}pB_{p-3}\pmod{p^2},
 \end{equation}
 and
 \begin{equation}\label{1.4}\sum_{k=0}^{\frac{p-3}2}\f{\bi{2k}k^2H_{2k}}{(2k+1)16^k}\eq-2\big(\frac{-1}p\big)E_{p-3}\pmod{p}.
 \end{equation}
\end{theorem}

  \begin{remark}\label{Rem1.1}
Our approach to Theorem 1.1 is somewhat unique in the sense that it depends heavily on some special combinatorial identities. And the reason why we research these two congruences is to prove the following Theorem. And the second congruence was conjectured by Z.-W. Sun in 2012, he also showed that it was equivalent to $$\sum_{k=1}^{(p-1)/2}\frac{\binom{2k}k^2}{k{16}^k}H_{2k}\equiv4\big(\frac{-1}p\big)E_{p-3}\pmod{p}.$$ He also list this conjecture in \cite{Su2}.
\end{remark}

  \begin{theorem}\label{Th1.2} For any odd prime $p$. We have
\begin{equation}\label{1.5}
\sum_{k=0}^{p-1}\frac{P_k}{8^k}\equiv1+2\big(\frac{-1}p\big)p^2E_{p-3}\pmod{p^3},
 \end{equation}
 and \begin{equation}\label{1.6}
 \sum_{k=0}^{p-1}\frac{P_k}{{16}^k}\equiv\big(\frac{-1}p\big)-p^2E_{p-3}\pmod{p^3}.
 \end{equation}
  \end{theorem}
 We are going to prove Theorem \ref{Th1.1} in section 2, and at last we prove Theorem \ref{Th1.2} in section 3.

  \section{Proof of Theorem 1.1}
 \setcounter{lemma}{0}
\setcounter{theorem}{0}
\setcounter{corollary}{0}
\setcounter{remark}{0}
\setcounter{equation}{0}
\setcounter{conjecture}{0}
  \begin{lemma}\label{Lem2.1}\cite[Lemma 4.2]{Su3} Let $p=2n+1$ be an odd prime, and let $k\in\{0,\ldots,n\}$. Then
  \begin{equation}\label{2.1}
  \frac{\binom nk}{\binom{2k}k/(-4)^k}\equiv1-p\sum_{j=1}^k\frac1{2j-1}\pmod{p^2}
  \end{equation}
  and
  \begin{equation}\label{2.2}
  \binom nk\binom{n+k}k(-1)^k\equiv\frac{\binom{2k}k^2}{16^k}\pmod{p^2}.
  \end{equation}
  \end{lemma}

   For any positive integer $n$, we have the following identities
   \begin{align}\label{2.3}
  \sum_{k=0}^{n-1}\frac{\binom nk\binom{n+k}k(-1)^k}{2k+1}H_k=-\frac{H_n\binom{2n}n(-1)^n}{2n+1}+\frac{2}{2n+1}\sum_{k=1}^n\frac{(-1)^k}k
  \end{align}
  and
  \begin{align}\label{2.4}
  \sum_{k=0}^{n-1}\frac{\binom nk\binom{2k}k}{(2k+1)(-4)^k}=-\frac{\binom{2n}n}{(2n+1)(-4)^n}+\frac{4^n}{(2n+1)\binom{2n}n}\sum_{k=0}^n\frac{\binom{2k}k^2}{16^k}
  \end{align}
These two identities can be easily proved by induction.
\begin{lemma}\label{Lem2.2}\cite[Theorem3.2]{S1}
\begin{align*}
\sum_{1\leq k<\frac p4}\frac 1k\equiv&-3q_p(2)+p(\frac32q_p(2)^2+(-1)^{(p-1)/2}(E_{2p-4}-2E_{p-3}))
\\&-p^2(q_p(2)^3+\frac7{12}B_{p-3})\pmod{p^3}
\end{align*}
and
\begin{align*}
\sum_{\frac p4<k<\frac p2}\equiv& q_p(2)-p(\frac12q_p(2)^2+(-1)^{(p-1)/2}(E_{2p-4}-2E_{p-3}))
\\&+\frac13p^2q_p(2)^3\pmod{p^3}.
\end{align*}
\end{lemma}

\begin{lemma}\label{Lem2.3} Let $p>3$ be a prime. Then
$$\binom{p-1}{(p-1)/2}\equiv(-1)^{(p-1)/2}4^{p-1}\pmod{p^3}.$$
\end{lemma}
\begin{remark}\label{Rem2.1}  Lemma \ref{Lem2.3} is a famous congruence of Morley\cite{Mo}.
\end{remark}
{\it Proof of theorem \ref{Th1.1}}: For any prime $p>3$, taking $n=(p-1)/2$ in (\ref{2.3}), by Lemma \ref{Lem2.2} and Lemma \ref{Lem2.3} we can get that
\begin{align*}
&\sum_{k=0}^{\frac{p-3}2}\frac{\binom{\frac{p-1}2}k\binom{\frac{p-1}2+k}k(-1)^kH_k}{2k+1}=-\frac{H_{\frac{p-1}2}\binom{p-1}{\frac{p-1}2}(-1)^{\frac{p-1}2}}p+\frac2p\sum_{k=1}^{\frac{p-1}2}\frac{(-1)^k}k
\\\eq&-\frac{H_{\frac{p-1}2}\binom{p-1}{\frac{p-1}2}(-1)^{\frac{p-1}2}}p+\frac2p(\sum_{1\leq k<\frac p4}\frac1k-\sum_{1\leq k<\frac p2}\frac1k)
\\\eq&\frac{(2q_p(2)-pq_p(2)^2+\frac23p^2q_p(2)^3+\frac7{12}p^2B_{p-3})(1+2pq_p(2)+p^2q_p(2)^2)}p
\\&+\frac2p(-q_p(2)+p(\frac12q_p(2)^2+(-1)^{\frac{p-1}2}(E_{2p-4}-2E_{p-3}))-\frac13p^2q_p(2)^3)
\\\eq&4q_p(2)^2+2(-1)^{\frac{p-1}2}(E_{2p-4}-2E_{p-3})+\frac7{12}pB_{p-3}\pmod{p^2}
\end{align*}
with (\ref{2.2}) we know $$\sum_{k=0}^{\frac{p-3}2}\frac{\binom{\frac{p-1}2}k\binom{\frac{p-1}2+k}k(-1)^kH_k}{2k+1}\equiv\sum_{k=0}^{\frac{p-3}2}\frac{\binom{2k}k^2 H_k}{(2k+1)16^k}\pmod{p^2},$$
so we finish the proof of (\ref{1.3}).

For any prime $p>3$ , taking $n=\frac{p-1}2$ in (\ref{2.4}), with Lemma \ref{Lem2.3} and (\ref{1.2}) we can deduce that
\begin{align*}
&\sum_{k=0}^{\frac{p-3}2}\frac{\binom{\frac{p-1}2}k\binom{2k}k}{(2k+1)(-4)^k}=-\frac{\binom{p-1}{\frac{p-1}2}}{p(-4)^{\frac{p-1}2}}+\frac{2^{p-1}}{p\binom{p-1}{\frac{p-1}2}}\sum_{k=0}^{\frac{p-1}2}\frac{\binom{2k}k^2}{16^k}
\\\eq&-\frac{2^{p-1}}p+\frac{(-1)^{\frac{p-1}2}}{p2^{p-1}}((-1)^{\frac{p-1}2}+p^2E_{p-3})
\\\eq&-\frac{2^{p-1}}p+\frac1{p2^{p-1}}+(-1)^{\frac{p-1}2}pE_{p-3}
\\\eq&-2q_p(2)+pq_p(2)^2+(-1)^{\frac{p-1}2}pE_{p-3}\pmod{p^2}
\end{align*}
(Note that $2^{p-1}=1+pq_p(2)$ and $4^{p-1}=(2^{p-1}-1+1)^2=(1+pq_p(2))^2=1+2pq_p(2)+p^2q_p(2)^2$.)

So we have
\begin{align}\label{2.5}
\sum_{k=0}^{\frac{p-3}2}\frac{\binom{\frac{p-1}2}k\binom{2k}k}{(2k+1)(-4)^k}\equiv-2q_p(2)+pq_p(2)^2+(-1)^{\frac{p-1}2}pE_{p-3}\pmod{p^2}.
\end{align}
Let $n=\frac{p-1}2$, it is known that $$\sum_{k=0}^{n-1}\frac{\binom{2k}k^2H_{2k}}{(2k+1)16^k}=\frac12\sum_{k=0}^{n-1}\frac{\binom{2k}k^2H_{k}}{(2k+1)16^k}+\sum_{k=0}^{n-1}\frac{\binom{2k}k^2}{(2k+1)16^k}\sum_{j=1}^k\frac1{2j-1}.$$
With (\ref{2.1}) we have
\begin{align*}
&p\sum_{k=0}^{n-1}\frac{\binom{2k}k^2}{(2k+1)16^k}\sum_{j=1}^k\frac1{2j-1}
\\\equiv&\sum_{k=0}^{n-1}\frac{\binom{2k}k^2}{(2k+1)16^k}-\sum_{k=0}^{n-1}\frac{\binom nk\binom{2k}k}{(2k+1)(-4)^k}\pmod{p^2}
\end{align*}
then
\begin{align*}
&p\sum_{k=0}^{n-1}\frac{\binom{2k}k^2H_{2k}}{(2k+1)16^k}
\\\equiv&\frac p2\sum_{k=0}^{n-1}\frac{\binom{2k}k^2H_{k}}{(2k+1)16^k}+\sum_{k=0}^{n-1}\frac{\binom{2k}k^2}{(2k+1)16^k}-\sum_{k=0}^{n-1}\frac{\binom nk\binom{2k}k}{(2k+1)(-4)^k}\pmod{p^2}
\end{align*}

by (\ref{1.1}), (\ref{1.3}) and (\ref{2.5}) we can deduce that
\begin{align*}
&p\sum_{k=0}^{n-1}\frac{\binom{2k}k^2H_{2k}}{(2k+1)16^k}
\\\eq&\frac p2(4q_p(2)^2-2(-1)^{\frac{p-1}2}E_{p-3})+(-2q_p(2)-pq_p(2)^2)
\\&-(-2q_p(2)+pq_p(2)^2+(-1)^{\frac{p-1}2}pE_{p-3})
\\\eq&-2(-1)^{\frac{p-1}2}pE_{p-3}\pmod{p^2}
\end{align*}
(Note that $E_{2p-4}\equiv E_{p-3}\pmod p$, \cite[(3.1)]{S1}).

Hence
$$\sum_{k=0}^{n-1}\frac{\binom{2k}k^2H_{2k}}{(2k+1)16^k}\equiv-2(-1)^{\frac{p-1}2}E_{p-3}\pmod p.$$
So we have done the proof of (\ref{1.4}).
 \section{Proof of Theorem 1.2}
 \setcounter{lemma}{0}
\setcounter{theorem}{0}
\setcounter{corollary}{0}
\setcounter{remark}{0}
\setcounter{equation}{0}
\setcounter{conjecture}{0}
{\it Proof of (\ref{1.5})}: For any prime $p>3$,  with Lemma \ref{Lem2.3} we have
\begin{align*}
&\sum_{k=0}^{p-1}\frac{P_k}{8^k}=\sum_{k=0}^{p-1}\sum_{j=0}^{\lfloor\frac k2\rfloor}\frac{\binom{k}{2j}\binom{2j}j^2}{16^j}
\\=&\sum_{j=0}^{\frac{p-1}2}\frac{\binom{2j}j^2}{16^j}\sum_{k=2j}^{p-1}\binom{k}{2j}=\sum_{j=0}^{\frac{p-1}2}\frac{\binom{2j}j^2}{16^j}\binom p{2j+1}
\\\eq&\frac{\binom{p-1}{\frac{p-1}2}^2}{4^{p-1}}+p\sum_{j=0}^{\frac{p-3}2}\frac{\binom{2j}j^2}{(2j+1)16^j}(1-pH_{2j})
\\\eq&4^{p-1}+p\sum_{j=0}^{\frac{p-3}2}\frac{\binom{2j}j^2}{(2j+1)16^j}-p^2\sum_{j=0}^{\frac{p-3}2}\frac{\binom{2j}j^2H_{2j}}{(2j+1)16^j}\pmod{p^3}
\end{align*}
(Note that $\binom{p-1}{2j}=\Pi_{l=1}^{2j}\frac{p-l}l\equiv1-pH_{2j}\pmod{p^2}$ for all $0\leq j\leq (p-3)/2$.)

So by (\ref{1.1}) we have
\begin{align*}
&\sum_{k=0}^{p-1}\frac{P_k}{8^k}\eq1+2pq_p(2)+p^2q_p(2)^2+p(-2q_p(2)-pq_p(2))
\\&-p^2\sum_{j=0}^{\frac{p-3}2}\frac{\binom{2j}j^2H_{2j}}{(2j+1)16^j}\eq1-p^2\sum_{j=0}^{\frac{p-3}2}\frac{\binom{2j}j^2H_{2j}}{(2j+1)16^j}\pmod{p^3},
\end{align*}
 then with (\ref{1.4}) we can easily get
$$\sum_{k=0}^{p-1}\frac{P_k}{8^k}\equiv1+2\big(\frac{-1}p\big)p^2E_{p-3}\pmod{p^3}.$$
 (\ref{1.5}) in the case of $p=3$ can be verified directly.
So we complete the proof of (\ref{1.5}).

{\it Proof of (\ref{1.6})}. First we know that in \cite{Su4} Z.-W.Sun proved that $$\sum_{k=0}^{p-1}\frac{\binom{2k}k^2}{16^k}\equiv\big(\frac{-1}p\big)-p^2E_{p-3}\pmod{p^3}.$$By taking $m=n$ and $x=\frac12$ in \cite[(1.1)]{TM} we can deduce that $$\sum_{k=0}^n\binom{n+k}k\frac1{2^k}=2^n.$$
With those two results we can get that
\begin{align*}
&\sum_{k=0}^{p-1}\frac{P_k}{{16}^k}=\sum_{k=0}^{p-1}\frac{S_k}{8^k}=\sum_{n=0}^{p-1}\frac1{8^n}\sum_{k=0}^n\binom{2k}k^2\binom{k}{n-k}(-4)^{n-k}
\\=&\sum_{n=0}^{p-1}\frac1{(-2)^n}\sum_{k=0}^n\binom{2k}k^2\binom k{n-k}\frac1{(-4)^k}=
\sum_{k=0}^{p-1}\frac{\binom{2k}k^2}{(-4)^k}\sum_{n=k}^{p-1}\frac{\binom{k}{n-k}}{(-2)^n}
\\=&\sum_{k=0}^{p-1}\frac{\binom{2k}k^2}{(-4)^k}\sum_{n=0}^{p-1-k}\frac{\binom kn}{(-2)^{n+k}}=
\sum_{k=0}^{p-1}\frac{\binom{2k}k^2}{8^k}\sum_{n=0}^{p-1-k}\frac{\binom kn}{(-2)^n}
\\=&\sum_{k=0}^{(p-1)/2}\frac{\binom{2k}k^2}{8^k}\sum_{n=0}^k\binom kn(-\frac12)^n+
\sum_{k=(p+1)/2}^{p-1}\frac{\binom{2k}k^2}{8^k}\sum_{n=0}^{p-1-k}\binom kn(-\frac12)^n
\\=&\sum_{k=0}^{(p-1)/2}\frac{\binom{2k}k^2}{8^k}(\frac12)^k+
\sum_{k=(p+1)/2}^{p-1}\frac{\binom{2k}k^2}{8^k}\sum_{n=0}^{p-1-k}\binom kn(-\frac12)^n
\end{align*}
 while for each $k=(p+1)/2,\ldots,p-1$, we have $\binom{2k}k^2\equiv0\pmod{p^2}$  and
 \begin{align*}
 &\sum_{n=0}^{p-1-k}\binom kn(-\frac12)^n=\sum_{n=0}^{p-1-k}\binom{-k+n-1}n(\frac12)^n
 \\&\equiv\sum_{k=0}^{p-1-k}\binom{p-1-k+n}n(\frac12)^n
=2^{p-1-k}\equiv\frac1{2^k}\pmod p
 \end{align*}
 so we have $$\sum_{k=(p+1)/2}^{p-1}\frac{\binom{2k}k^2}{8^k}\sum_{n=0}^{p-1-k}\binom kn(-\frac12)^n\equiv\sum_{k=(p+1)/2}^{p-1}\frac{\binom{2k}k^2}{16^k}\pmod{p^3}.$$
 Hence $$\sum_{k=0}^{p-1}\frac{P_k}{{16}^k}\equiv\sum_{k=0}^{p-1}\frac{\binom{2k}k^2}{16^k}\equiv\big(\frac{-1}p\big)-p^2E_{p-3}\pmod{p^3}.$$
 Thus we complete the proof of (\ref{1.6}).

 Until now we finish the proof of Theorem \ref{Th1.2}.
 
\end{document}